\documentclass[12pt,reqno]{amsart}
\usepackage{graphicx} 
\usepackage{tikz,xypic}
\usepackage{url, hyperref}
\usepackage{fullpage}
\usepackage{microtype}
\usepackage{amsmath, amssymb, amsthm}
\usepackage{listings}
\usepackage{color}

\definecolor{mygreen}{rgb}{0,0.6,0}
\definecolor{mygray}{rgb}{0.5,0.5,0.5}
\definecolor{mymauve}{rgb}{0.58,0,0.82}
\definecolor{backcolour}{rgb}{0.95,0.95,0.92}
\theoremstyle{definition}
    \newtheorem{theorem}{Theorem}
    
    \newtheorem{lemma}[theorem]{Lemma}
    \newtheorem{proposition}[theorem]{Proposition}
    \newtheorem{remark}[theorem]{Remark}
    
    \newtheorem{definition}[theorem]{Definition}
    \newtheorem{conjecture}[theorem]{Conjecture}

\newcommand{\bfnu}{\boldsymbol{\nu}}
\newcommand{\bfmu}{\boldsymbol{\mu}}

\DeclareMathOperator{\Aut}{Aut}

\DeclareMathOperator{\val}{val}

\DeclareMathOperator{\mult}{mult}
\DeclareMathOperator{\trop}{trop}
\def\bfx{\mathbf{x}}
\def\bfe{\mathbf{e}}
\DeclareMathOperator{\DR}{DR}
\DeclareMathOperator{\branch}{branch}

\title{One part leaky covers}
\author{Renzo Cavalieri} 
\address{Colorado State University, Department of Mathematics, Weber Building, Fort Collins, CO 80523-1874, USA}
\email{renzo@math.colostate.edu}

\author {Hannah Markwig}
\address {Universit\"at T\"ubingen, Fachbereich Mathematik, Auf der Morgenstelle 10, 72076 T\"ubingen, Germany }
\email {hannah@math.uni-tuebingen.de}

\author{Johannes Schmitt}
\address{Department Mathematik, R\"amistrasse 101, CH-8092 Z\"urich, Switzerland}
\email{johannes.schmitt@math.ethz.ch}

\DeclareMathOperator{\logDR}{logDR}

\newcommand{\Mgn}{\overline{\mathcal{M}}_{g,n}}
\begin{document}

\maketitle

\begin{abstract}
   In our previous work \cite{CMS24} we defined a new class of enumerative invariants called  $k$-leaky double Hurwitz descendants, generalizing both descendant integrals of double ramification cycles and $k$-leaky double Hurwitz numbers. Here, we focus on the one-part version of these numbers, i.e.\ when the positive ramification profile is $(d)$. We derive recursions and use them to produce explicit formulas and structure results for some infinite families of these numbers.
\end{abstract}



\tableofcontents

\section{Introduction}

 We begin this introduction by summarizing the original results of this work, which we then comment on in the second section. In the third section we expand our discussion to provide an overview of how the problem fits within its broader field of research. While this ``onion-style'' approach is non-standard, we feel it provides the most efficient access to the contents of this paper for readers that are somewhat familiar with the interactions of tropical and logarithmic geometry in curve counting problems. We tried to write each section in a self-contained way, so readers may read them in whatever order suits them best.

\subsection{Results}

We study one-part $k$-leaky double Hurwitz numbers with $\psi$-class insertions $\mathsf{H}_g((d, - \bfnu), \mathbf{e})$, with the main goal of exhibiting explicit formulas and unveiling part of the rich combinatorial structure of these enumerative invariants. We will  give precise definitions later, but for now it suffices to view $\mathsf{H}_g((d, - \bfnu), \mathbf{e})$ as an enumerative invariant for pointed curves of genus $g$; the discrete data $(d, - \bfnu)$ and $\mathbf{e}$ encodes ramification conditions and  $\psi$-class insertions at the marked points.


We first consider leaky descendant Hurwitz numbers with genus $g=0$ and $\psi$-insertions only at the unique positive marking $p_0$ and at most one other marking, which we choose to be $p_1$.  We introduce the notation
\begin{equation} \label{eqn:def_J_e0e1}
     J_{e_0,e_1}^m  = J_{e_0,e_1}^m (k,d,\nu_1):=\mathsf{H}_0((d, - \bfnu), (e_0, e_1, 0,\ldots, 0))
\end{equation}
As suggested by the notation, these invariants are independent from the ramification orders at points that do not support a $\psi$-class. 

\begin{theorem} \label{Thm:genus0}
The following formulas hold for the families of $k$-leaky one-part double Hurwitz descendants  defined in \eqref{eqn:def_J_e0e1}:
\begin{align}
    J^m_{e_0,0}(k,d,\nu_1) &= \frac{(m-1)!}{(e_0+1)!} \cdot \prod_{j=e_0+1}^{m-2} \left( d-j\cdot\frac{k}{2}\right)\,, \label{eqn:J_e0_0}\\
    \intertext{ and}
    J^m_{0,e_1}(k,d,\nu_1) &= \frac{(m-2)!}{(e_1+1)!} \cdot \Big( e_1 \cdot \frac{\prod_{j=e_1+1}^{m-1} (d-\nu_1-j\cdot\frac{k}{2}) - \prod_{j=1}^{m-e_1-1} (d-j\cdot\frac{k}{2})}{-\nu_1-e_1\cdot\frac{k}{2}} \notag\\ 
    &\quad + (m-e_1-1) \cdot \prod_{j=1}^{m-e_1-2} \left(d-j\cdot\frac{k}{2}\right) \Big)\,\label{eqn:J_0_e1}.
\end{align}

\end{theorem}
Formula \eqref{eqn:J_0_e1}  has the following specializations  when setting $k=0$ or $\nu_1=0$:
\begin{align}
    J^m_{0,e_1}(0,d,\nu_1)  &= \frac{(m-2)!}{(e_1+1)!} \cdot \left( e_1 \cdot  \frac{d^{m-1-e_1} - (d-\nu_1)^{m-1-e_1}}{\nu_1}+(m-1-e_1)\cdot d^{m-2-e_1}\right)\,,\label{eqn:J_0_e1_k0}\\
      J^m_{0,e_1}(k,d,0)  &= \frac{(m-2)!}{(e_1+1)!} \cdot (m-1-e_1) \cdot \left( \sum_{i=1}^{e_1+1}  \ \prod_{j=0}^{m-e_1-3}\left(d-(i+j)\cdot\frac{k}{2}\right)\right)\,. \label{eqn:J_0_e1_nu10}
\end{align}

In Proposition \ref{Prop:onepart_genus0_structure} we prove, still for $g=0$, the structural result that for arbitrary $\textbf{e}$, the invariant $\mathsf{H}_0((d, - \bfnu), \mathbf{e})$ is a polynomial of degree $m-2-|\textbf{e}|$ in $d, k$ and  variables $\nu_i$ associated to markings with non-trivial $\psi$-insertions (i.e. indices $1 \leq i \leq m$ with $e_i > 0$).

In genus one, we prove a formula that describes one part leaky Hurwitz numbers (with no descendant insertions) for arbitrary values of $k$.

\begin{theorem}\label{thm:genusone}
The following formula holds:

\begin{align} \label{eqn:GJV_genus_1}
    & \mathsf{H}_1((d, - \bfnu), \textbf{0})=\\& \nonumber\frac{(m+1)!}{24} \cdot  \prod_{j=2}^m \left(d - j\cdot \frac{k}{2}\right)
    \cdot\left[
(d-k)^3 - \left(d-\frac{k}{2}\right)\cdot\left(1+ 2\cdot \sum_{s<t}(\nu_s+k)(\nu_t+k)\right)\right].   
\end{align}
  \end{theorem}

In arbitrary genus $g \geq 0$ we have a complete understanding of the $k=0$ specialization of the $k$-leaky double Hurwitz descendants with all $\psi$-insertions at marking $p_0$; the structure we observe naturally generalizes and interpolates between the case $e_0=0$ treated in \cite[Theorem 3.1]{GJV05} and the case $e_0=2g-2+m$ from \cite[Theorem 1]{BSSZ15}.

\begin{theorem}  \label{Thm:psi0_e_general_g}
For $g,m \geq 0$ with $2g-1+m>0$, 
we have
\begin{equation}
    \mathsf{H}_g((d,-\bfnu), (e_0, 0, \ldots, 0))|_{k=0} = \frac{d^{-e_0}}{(e_0+1)!}\cdot \mathsf{H}_g((d,-\bfnu), ( 0, \ldots, 0))|_{k=0}.
\end{equation}
\end{theorem}

In particular, since \cite{BSSZ15} provides a formula for  descendant invariants, we also obtain explicit formulas for all intermediate invariants.
Denote
\[
S(z) = \frac{\sinh(z/2)}{z/2} = \sum_{k \geq 0} \frac{z^{2k}}{2^{2k}(2k+1)!} = 1 + \frac{z^2}{24} + \frac{z^4}{1920} + \frac{z^6}{322560} + \ldots
\]
considered as a formal power series in $\mathbb{Q}[z]$.
We have
\begin{equation}
    \mathsf{H}_g((d,-\bfnu), (e_0,0,\ldots, 0))|_{k=0} = \frac{(2g-1+m)!}{(e_0+1)!} \cdot d^{2g-2+m-e_0} \cdot [z^{2g}] \frac{\prod_{j=1}^m S(\nu_j z)}{S(z)}\,,
\end{equation}

where $[z^{2g}]$ means taking the coefficient of the monomial $z^{2g}$ in the  power-series expansion in $z$ of the expression that follows.

%

\subsection{Discussion}

Leaky descendant Hurwitz numbers, defined in \cite[Definition 1.1]{CMS24} are degrees of $0$-dimensional classes in the logarithmic Chow ring of $\Mgn$, obtained by multiplying a log DR cycle with two distinguished types of classes: $\psi$-classes pulled-back from $\Mgn$ and piecewise polynomial branch classes pulled back from  $\mathsf{tEx}$. Consider $\mathbf{x} = (d, - \bfnu) =  (d, -\nu_1, \ldots, -\nu_m)$ of length $n=m+1$, where $d, \nu_1, \ldots, \nu_m > 0$ satisfy
$$d - \sum_{i=1}^m \nu_i = k \cdot(2g-2+n)\,.$$Labeling  the marks of $\overline{\mathcal{M}}_{g,n}$ by $p_0, p_1, \ldots, p_m$ and denoting $e = |\mathbf{e}|$,   we have
\begin{equation}\label{eq:leakydesdef}
\mathsf{H}_g((d, - \bfnu), \mathbf{e}) = \int_{\logDR_g(d, - \bfnu)} \psi_0^{e_0}\cdots \psi_m^{e_m} \cdot \branch_{2g-3+n-e}\,.
\end{equation}

This paper takes the first steps in studying the algebraic and combinatorial structure arising from the dependence of these enumerative invariants on their discrete data. We deem this endeavor interesting for two reasons: first because, when $k = 0$, these invariants naturally interpolate between double Hurwitz numbers and descendant invariants of the DR cycle, two classical families of enumerative invariants that are well-known to have rich structure. Secondly, letting $k>0$ one can observe interesting parallels between invariants of moduli spaces of maps and of moduli spaces of (pluri)differentials.

The main tool we use to study this structure is a tropical algorithm (\cite[Theorem 1.2]{CMS24}) computing leaky descendant Hurwitz numbers as a weighted sum over leaky tropical covers, see Section \ref{sec:background}.

Such perspective provides a natural recursion among these invariants: truncating all leaky covers before their last vertex and organizing the sums in terms of the type of the last vertex expresses a leaky descendant Hurwitz number as a finite sum of similar invariants with smaller value for the quantity $2g-2+n$. This recursion is stated in Lemmas \ref{Lem:genus_0_recursion0}, \ref{Lem:genus_0_recursion}, \ref{lem:reggenzerogenerale} and \ref{prop-rechighergenus}, each time in the generality needed at that point.

While such recursive approach is always available, rather than attempting a general description, we choose to focus our attention on a few infinite families of double leaky Hurwitz descendants where the structure is especially explicit and simple.

Restricting to the one-part chamber (where we have only the point $p_0$ with positive ramification condition) has two simplifying effects:  the recursion only involves connected invariants, and the resulting families of invariants are polynomial, rather than piecewise polynomial.

Limiting the number of $\psi$-insertions has the effect of limiting the combinatorial possibilities for the cut-off vertex in the tropical leaky covers, which reduces the combinatorial complexity of the recursion.

In genus zero, the last vertex of any tropical leaky cover must have a unique left pointing flag, so all recursive terms have strictly fewer marked points. We focus our attention on invariants with exactly one descendant insertion of arbitrary power. It is convenient to introduce the notation $J_{e_0,e_1}^m$  \eqref{eqn:def_J_e0e1} that allows insertions both on the unique mark with positive profile as well as on a single mark with negative profile in order to state efficiently the recursion; eventually Theorem \ref{Thm:genus0} gives explicit formulas for these invariants in the cases $(e_0, 0)$ and $(0,e_1)$.

Our path towards Theorem \ref{Thm:genus0} went through the following four steps:
\begin{enumerate}
    \item
    The tropical algorithm computing the numbers $ J_{e_0,e_1}^m $ was implemented based on the software package \texttt{admcycles} \cite{admcycles}. From \cite[Theorem 1.3]{CMS24} it follows that $ J_{e_0,e_1}^m $ is a polynomial of degree $m-2-e_0-e_1$ in $d, \nu_1, \ldots, \nu_n$, where $k$ is then determined as 
    $$k=(d-\nu_1 - \ldots - \nu_m)/(m-1).$$
    Thus for any fixed $e_0, e_1, m$, a formula for $ J_{e_0,e_1}^m $ can be obtained using multivariate interpolation by evaluating at finitely many points $\mathbf{x} = (d, - \bfnu)$.
    \item Encouraged by the simple shape of the resulting formulas (for small parameters $e_0, e_1, m$), we identified the  recursion for $ J_{e_0,e_1}^m $, see Lemma \ref{Lem:genus_0_recursion}. This allows the calculation for larger values of $m$, which would be infeasible to obtain just from the interpolation method of Step 1.
    \item In the cases $e_1=0$ or $e_0=0$ presented in Theorem \ref{Thm:genus0}, we  guessed the shape of the formulas presented there from the computational evidence and some intuition derived from the shape of the recursion in Lemma \ref{Lem:genus_0_recursion}. In the case $e_0=0$, we first found the specializations \eqref{eqn:J_0_e1_k0} and \eqref{eqn:J_0_e1_nu10} before managing to find the general shape \eqref{eqn:J_0_e1}.
    \item Once a candidate formula for $ J_{e_0,e_1}^m $ is available, it can be proven to satisfy the recursion derived in Lemma \ref{Lem:genus_0_recursion} via lengthy, but ultimately manageable algebraic manipulations.
\end{enumerate}
In generalizing the process described above, the bottleneck (for the authors) lies in Step 3, where a general formula must be guessed from limited data. This data is available in many further cases not covered above (e.g. when both $e_0, e_1 > 0$), and we expect that reasonably simple formulas could exist here. 

In genus one we restrict our attention to the case $\mathbf{e} =0$, where the formula we guessed is shown  to satisfy the recursion from Lemma \ref{prop-rechighergenus} by expressing it as a polynomial in the elementary symmetric functions in the quantities $(\nu_i+k)$. The formula in Theorem \ref{thm:genusone} is easily shown to be equivalent to
\begin{align*}
   & \mathsf{H}_1((d, - \bfnu), \textbf{0}) =\\& \frac{(m+1)!}{24} \cdot \left( \left( \sum_{i=1}^{m} (\nu_i + k)^2 - 1 \right) \cdot \prod_{j=1}^{m} \left(d - j\cdot\frac{ k}{2}\right)  - k\cdot\frac{\left(d - k\right)^2}{2} \cdot  \prod_{j=2}^{m} \left(d - j\cdot\frac{ k}{2}\right)  \right)\,.
\end{align*}
which, setting $k=0$, immediately reduces to the formula for one part genus one double Hurwitz numbers
proven in \cite{GJV05}.

For arbitrary genus $g$, we set $k=0$ but allow one descendant insertion of arbitrary degree at the unique point of positive ramification. In this case we have a different recursion at our disposal, consisting of exchanging one $\psi$-class for a linearly equivalent boundary divisor. In $k = 0$, such operation follows from the existence of a branch morphism from the DR cycle to a (quotient of) Losev-Manin space, as well as known comparison results among the $\psi$-classes on Losev-Manin, DR and $\Mgn$. In \cite[Proposition 4.2]{CMS24} we  show a similar exchanging rule that trades a $\psi$-class with a boundary divisor plus a multiple of the class $\kappa_1$ restricted to the logarithmic double ramification cycle.  

In Section \ref{Sect:Implementations} we describe the existing computer implementations of formulas for $ \mathsf{H}_g(\mathbf{x}, \mathbf{e})$ (and how to access them). The authors  encourage interested readers to try their hands at finding new formulas and shedding light on the intriguing structure that can be already observed in these first steps. For example, there seems to be a hypergeometric dependence of these invariants on the parameter $k$, and it would be satisfying to have a conceptual explanation for this phenomenon.

\subsection{Context}
Hurwitz numbers are classical enumerative invariants counting branched covers of the projective line \cite{hurwitz}. The ELSV formula \cite{ELSV01}  establishes a  connection between single Hurwitz numbers and intersection theory on moduli spaces of curves, revealing polynomial structure and links to Gromov–Witten theory.

Double Hurwitz numbers $H_g(\bfmu,\bfnu)$ extend this setting by counting genus $g$ covers with specified ramification profiles $\bfmu$ over $0$ and $\bfnu$ over $\infty$, with simple branching elsewhere. The seminal work of Goulden, Jackson, and Vakil \cite{GJV05} established that for fixed genus $g$, these numbers are piecewise polynomial in the parts of the partitions $\bfmu$ and $\bfnu$. Shadrin, Shapiro, and Vainshtein \cite{SSV08} analyzed the chamber structure of $H_0(\bfmu,\bfnu)$, explicitly describing wall-crossing behavior, while Cavalieri, Johnson, and Markwig \cite{CJM11} developed tropical methods that provided a systematic derivation of wall-crossing formulas in all genera. Johnson \cite{Joh15} employed infinite-wedge Fock space formalism to derive explicit formulas for double Hurwitz numbers with arbitrary ramification profiles, providing another proof of piecewise polynomiality.


A  special chamber of polynomiality  parameterizes one-part double Hurwitz numbers, where one ramification profile consists of a single part.
In genus 0, Goulden-Jackson-Vakil established the  formula:
\begin{equation}
{H_0}((d),\bfnu) = (n-1)! \cdot d^{n-2}
\end{equation}
where $\bfnu=(a_1,\ldots,a_n)$ is a partition of degree $d$ with $n$ parts. For genus $g\leq 5$, they compute explicit polynomial expressions representing  one part double Hurwitz numbers as the product of a monomial in $d$ and even symmetric polynomials in the $\nu_i$'s.
They further conjectured an ELSV-type formula for one-part double Hurwitz numbers $H_g((d),\bfnu)$, i.e. a tautological intersection theoretic expression on some compactification of the universal Picard stack of the form:
\begin{equation}\label{eq:elsv}
 H_g((d),\bfnu)   = d\cdot r! \int_{\overline{\mathrm{Pic}}_{g,n}} \frac{\Lambda_0 -\Lambda_2 +\ldots +(-1)^g\Lambda_{2g}}{\prod 1 - \nu_i\psi_i},
\end{equation}
where $r$ denotes the number of simple branch points, the $\psi_i$'s should be closely related to the pull-back of the homonymous classes on $\overline{\mathcal{M}}_{g,n}$ and the $\Lambda_{2i}$ are some yet unidentified classes of degree $2i$. Expression \eqref{eq:elsv} renders various properties (degree, parity, and vanishing of coefficients in certain degree ranges)  of the one-part Hurwitz polynomial transparent.

While the Goulden-Jackson-Vakil conjecture in its strong form $\eqref{eq:elsv}$ is still wide open,  some results \cite{CM14,DoLe22}  expressed one-part double Hurwitz numbers as intersection numbers on moduli spaces of curves, confirming the predicted polynomial and integrability structure. Blot \cite{Blot22} established a connection between one-part double Hurwitz numbers and the quantum KdV hierarchy, showing that the coefficients of $H_g((d),\bfnu)$ coincide with the expansion of the quantum Witten–Kontsevich series.


Cavalieri, Markwig, and Ranganathan \cite{CMR22} introduced $k$-leaky double Hurwitz numbers $H_g^k(\bfmu,\bfnu)$, which relax the usual degree-balancing condition by allowing a controlled discrepancy: $\sum \mu_i - \sum \nu_j = k \cdot(2g-2)$. These invariants are defined via logarithmic geometry as intersection numbers against pluricanonical double ramification cycles. They generalize ordinary double Hurwitz numbers (recovered when $k=0$) and maintain key properties: for fixed $g$ and $k$, they are piecewise polynomial in the parts of $\bfmu$ and $\bfnu$, with well-defined chamber walls. 

Descendant insertions, or $\psi$-classes, are natural elements of the Chow ring of $\Mgn$. The structure of their intersection numbers was unveiled in Witten's conjecture (Kontsevich's theorem)\cite{witten}, stating that a generating function for all intersection numbers of $\psi$-classes is a $\tau$-function for the KdV hierarchy.
Buryak, Shardin, Spitz and Zvonkine \cite{BSSZ15} describe the structure of intersection numbers of $\psi$-classes against the double ramification cycle.

Cavalieri, Markwig, and Schmitt \cite{CMS24}  defined  $k$-leaky double Hurwitz descendants, incorporating $\psi$-classes at the special ramification points. These descendant invariants also exhibit piecewise polynomiality and satisfy explicit wall-crossing formulas in genus 0.


Recent work by Accadia, Karev, and Lewański \cite{AKL25} approaches $k$-leaky Hurwitz numbers from the bosonic Fock space perspective, extending the basis of symmetric functions to include completed cycle operators. This provides a new proof of piecewise polynomiality in the leaky setting and suggests a potential ELSV-type formula involving completed cycles and quantum intersection numbers.


\section*{Acknowledgements}

We would like to thank Dawei Chen, Reinier Kramer, Danilo Lewański and Dhruv Ranganathan for conversations related to this project.

R.C. was supported by the Simon's foundation through MPS-TSM-00007937.
H.M. acknowledges support by the Deutsche Forschungsgemeinschaft (DFG, German Research Foundation), Project-ID 286237555, TRR 195. 
J.S. was supported by SNF-200020-219369 and SwissMAP.

\section{Background: tropical leaky covers}
\label{sec:background}

To set the notation for tropical leaky covers, we follow the exposition of \cite{CMS24}:
{an \emph{abstract tropical} \emph{curve} is a connected metric graph $\Gamma$ whose leaf edges (called \emph{ends}) have infinite length, together with a genus function $g:\Gamma\rightarrow \mathbb{Z}_{\geq 0}$ with finite support. Locally around a point $p$, $\Gamma$ is homeomorphic to a star with $r$ half-edges. 
The number $r$ is called the \emph{valence} of the point $p$ and denoted by $\val(p)$. The \emph{minimal vertex set} of $\Gamma$ is defined to be the points where the genus function is non-zero, together with points of valence different from $2$. 
Besides \emph{edges}, we introduce the notion of \emph{flags} of $\Gamma$. A flag is a pair $(\mathrm{v,e})$ of a vertex $\mathrm{v}$ and an edge $\mathrm{e}$ incident to it ($\mathrm{v}\in \partial \mathrm{e}$). Edges that are not ends are required to have finite length and are referred to as \emph{bounded} or \textit{internal} edges.

A \emph{marked tropical curve} is a tropical curve whose leaves are labeled. An isomorphism of a tropical curve is an isometry respecting the leaf markings and the genus function. The \emph{genus} of a tropical curve $\Gamma$ is given by
\[
g(\Gamma) = h_1(\Gamma)+\sum_{p\in \Gamma} g(p).
\]
The \emph{combinatorial type} is the equivalence class of tropical curves obtained by identifying any two tropical curves which differ only by edge lengths.

We study covers of $\mathbb{R}$ by graphs up to additive translation,  and  equip  $\mathbb{R}$ with a polyhedral subdivision to ensure the result is a map of metric graphs (see e.g.\ Section 5.4 and \ Figure 3 in \cite{MW17}). A \textit{metric line graph} is any metric graph obtained from a polyhedral subdivision of $\mathbb{R}$. The metric line graph determines the polyhedral subdivision up to translation. We fix an orientation of a metric line graph going from left to right (i.e.\ from negative values in $\mathbb{R}$ to positive values).

\begin{definition}[Leaky cover, \cite{CMR22}] \label{Def:leakycov}
Let $\pi:\Gamma\rightarrow T$ be a surjective map of metric graphs where $T$ is a metric line graph. We require that $\pi$ is piecewise integer affine linear, the slope of $\pi$ on a flag or edge $\mathrm{e}$ is a positive integer called the \emph{expansion factor} $\omega(\mathrm{e})\in \mathbb{N}_{\geq 0}$. 

For a vertex $\mathrm{v}\in \Gamma$, the \emph{left (resp.\ right) degree of $\pi$ at $\mathrm{v}$} is defined as follows. Let $f_l$ be the flag of $\pi(\mathrm{v})$ in $T$ pointing to the left  ($f_r$ the flag pointing to the right). Add the expansion factors of all flags $f$ adjacent to $v$ that map to $f_l$ (resp.\ $f_r$):
\begin{equation}
d_\mathrm{v}^l=\sum_{f\mapsto f_l} \omega(f), \;\;\;\; d_\mathrm{v}^r=\sum_{f\mapsto f_r} \omega(f).
\end{equation} 

\noindent
We say that the $k$-leaky condition is satisfied at $\mathrm{v}\in \Gamma$ if
\begin{equation}\label{eq:leakybal}d_\mathrm{v}^l-d_\mathrm{v
}^r= k \cdot(2g(\mathrm{v})-2+\val(\mathrm{v})).\end{equation}

We impose a stability condition:  $\Gamma\to T$ is called stable if the preimage of every vertex of $T$ contains a vertex of $\Gamma$ in its preimage which is of genus greater than $0$ or valence greater than $2$. 

Furthermore, we stabilize the source tropical curve further by passing to its minimal vertex set (containing only the points where the genus function is non-zero, together with points of valence different from $2$). 
The outcome $\pi:\Gamma\rightarrow T$ is called a \emph{$k$-leaky cover}.
\end{definition}
By the stabilization procedure, we lose the property that the cover is a map of graphs, however, this vertex structure is relevant to determine valencies correctly for the purpose of $\psi$-conditions.

\begin{definition}[Left and right degree]
The \emph{left (resp.\ right) degree}  of a leaky cover is the tuple of expansion factors of its ends mapping asymptotically to $-\infty$ (resp.\ $+\infty$). 
The tuple is indexed by the labels of the ends mapping to  $-\infty$ (resp.\ $+\infty$).  When the order imposed by the labels of the ends plays no role, we drop the information and treat the left and right degree only as a multiset. 
\end{definition}
By convention, we denote the left degree by $\mathbf{x}^{+}$ and the right degree by $\mathbf{x}^{-}$. In the right degree, we use negative signs for the expansion factors, in the left degree positive signs. We also merge the two to one vector which we denote $\bfx=(x_1,\ldots,x_n)$  called the \emph{degree}. The labeling of the ends plays a role: the expansion factor of the end with the label $i$ is $x_i$. In $\mathbf{x}$, we  distinguish the expansion factors of the left ends from those of the right ends by their sign.

We focus on one part leaky covers, i.e. when $\mathbf{x}=(d, -\bfnu)$. In this case we let $n = m+1$ and index the components of $\bfnu$ from $1$ to $m$.
An Euler characteristic calculation, combined with the leaky cover condition \eqref{eq:leakybal}, shows that 
$$ d-\sum_{i=1}^m \nu_i=k\cdot(2g-2+n), $$
where $g$ denotes the genus of $\Gamma$. An automorphism of a leaky cover is an automorphism of $\Gamma$ compatible with $\pi$.}

For the remaining part of the section, let $g \geq 0$ such that $2g-2+n > 0$ and consider vectors $(d, -\bfnu) $ such that $d-\sum_{i=1}^m\nu_i = k \cdot (2g-2+n)$ for some $k \in \mathbb{Z}$. Let $\bfe \in \mathbb{Z}_{\geq 0}^n$ such that $0 \leq |\bfe| \leq 2g-3+n$.

\begin{definition}[$\psi$-conditions for leaky covers] \label{def:psi}
 Let $\pi:\Gamma\rightarrow T$ be a $k$-leaky cover.
For a vertex $\mathrm{v}$, let $I_\mathrm{v}\subset \{1,\ldots,n\}$ be the subset of ends adjacent to $\mathrm{v}$ after passing to the minimal vertex set of $\Gamma$.
We say that $\pi:\Gamma\rightarrow T$  {\it satisfies the $\psi$-conditions} $\mathbf{e}$
  if for all vertices $\mathrm{v}$ of $\Gamma$,
  \begin{equation} \label{eq:psicondition}
 \sum_{i\in I_\mathrm{v}}e_i =2g(\mathrm{v})-3+ \val(\mathrm{v}).\end{equation}
\end{definition}
\begin{definition}[Vertex multiplicities]\label{def-multvertex}
Let $\pi:\Gamma\rightarrow T$ be a $k$-leaky cover satisfying the $\psi$-conditions $\bfe$.
For a vertex $\mathrm{v}$, let $I_\mathrm{v}\subset \{1,\ldots,n\}$ be the subset of ends adjacent to $\mathrm{v}$ after passing to the minimal vertex set of $\Gamma$.
   Let $\mathbf{x}(\mathrm{v})$ denote the vector containing the (left and right) local degree of $\mathrm{v}$, let $g(\mathrm{v})$ denote the genus of $\mathrm{v}$.

   We define the \emph{vertex multiplicity} to be

   $$\mult_\mathrm{v}:= \int_{\overline{M}_{g(\mathrm{v}),\val(\mathrm{v})}} \DR_{g(\mathrm{v})}(\mathbf{x}(\mathrm{v}))\cdot \prod_{i\in I_\mathrm{v}} \psi_i^{e_i}$$
\end{definition}

\begin{definition}[Count of $k$-leaky covers satisfying $\psi$-conditions]\label{def-trophgxe}


We define
\begin{equation}
  \mathsf{H}^{\trop}_g((d, -\bfnu), \bfe) = \sum_\pi \mult_\pi,  
\end{equation}
where:

\begin{itemize}
    \item $\pi:\Gamma\rightarrow T$ ranges among all leaky covers of degree $\mathbf{x}$ and genus $g$ (Definition \ref{Def:leakycov}) and satisfying the $\psi$-conditions $\mathbf{e}$ (Definition \ref{def:psi}); we require that every vertex of $T$ has precisely one vertex in its preimage;
    \item the multiplicity
\begin{equation}
   \mult_\pi =  \frac{1}{|\Aut(\pi)|}\cdot \prod_\mathrm{e} \omega(\mathrm{e})\cdot \prod_\mathrm{v} \mult_\mathrm{v} \in \mathbb{Q}
\end{equation}
where the first part is the product of the expansion factors at the bounded edges of $\Gamma$ (according to its minimal vertex set), weighted by the number of automorphisms of $\pi$; the last  product goes over the set of vertices of $\Gamma$ with $\mult_\mathrm{v}$  as in Definition \ref{def-multvertex}.
\end{itemize}
\end{definition}

\begin{remark}[\cite{CMS24}, Remark 5.2]\label{rem-localvertexmultgenus0}
     In genus $0$, the vertex multiplicity 
    \[ \mult_\mathrm{v}= \frac{(\val(\mathrm{v})-3)!}{\prod_{i\in I_\mathrm{v}}e_i!}\] 
     of a leaky cover satisfying $\psi$-conditions does not depend on $k$ or on the expansion factors of its adjacent edges; it equals a multinomial coefficient that only depends on its valence and the $\psi$-conditions.
\end{remark}

\begin{theorem}[\cite{CMS24}, Theorem 1.2]\label{thm-corres}


The $k$-leaky double Hurwitz descendant \eqref{eq:leakydesdef} equals the count of tropical $k$-leaky covers satisfying $\psi$-conditions (Definition \ref{def-trophgxe}):

$$\mathsf{H}^{\trop}_g((d, -\bfnu), \bfe)=\mathsf{H}_g((d, -\bfnu), \bfe).$$

\end{theorem}

\section{Genus zero}
We study the case of rational invariants. Throughout this section, we set $g=0$ and
consider $(d, - \bfnu) =  (d, -\nu_1, \ldots, -\nu_m)$ of length $n=m+1$, where $d, \nu_1, \ldots, \nu_m > 0$ satisfy
\begin{equation}\label{eq:leakyglobbal}
d - \sum_{i=1}^m \nu_i = k\cdot(m-1)\,.\end{equation}

\subsection{Recursion and formulas}
We first focus on invariants with at most two descendant insertions. In this case we use the shortened notation
$$J_{e_0,e_1}^m=\mathsf{H}_0((d,-\nu_1,\ldots,-\nu_m),(e_0,e_1,0,\ldots,0)),$$
where we observe that  $J_{e_0,e_1}^m$ is a function of $d,\nu_i,k,e_0,e_1$ and $m$, but we leave visible only the variables that  participate in the recursions we use. The next two lemmas develop recursions among these invariants.

\begin{lemma} \label{Lem:genus_0_recursion0}
Given $m \geq 2$ and $0\leq e_0 \leq m-2$, we have
\begin{equation}\label{eq:simplifiedrec}
  J_{e_0,0}^m=(m-1)\cdot \left(d-(m-2)\cdot\frac{k}{2}\right)\cdot J_{e_0,0}^{m-1} +\delta_{m-2, e_0}.  
\end{equation} 
with  the base cases of the recursion encoded by the   Kronecker delta $\delta_{m-2,e_0}$.
In particular the invariants $J_{e_0, 0}^m$ are independent of $\bfnu$, the negative ramification profile.
\end{lemma}
\begin{proof}
We prove the independence of $\bfnu$ by induction on $m$, with the base case $m=2$ being true because the invariant $J^2_{0,0} = 1$. This also establishes \eqref{eq:simplifiedrec} for $m=2$.
For $m>2$, since $J^m_{m-2,0} = \int_{\overline{\mathcal{M}}_{0,m+1}}\psi_0^{m-2} = 1$, \eqref{eq:simplifiedrec} holds for $e_0 = m-2$.

 By Theorem \ref{thm-corres}, the numbers $J_{e_0, 0}^m$ are calculated as a weighted count of tropical $k$-leaky covers, which for $e_0<m-2$ have at least two vertices. To establish the recursive part of \eqref{eq:simplifiedrec},we cut each  tropical cover just before  the last (rightmost) vertex, which is $3$-valent and has multiplicity $1$. If it is adjacent to the ends of weight $\nu_i$ and $\nu_j$,  the weight of its adjacent bounded edge is $\nu_i+\nu_j+k$. Cutting this edge, the remainder of the graph (to the left of the cut) contributes to $J_{e_0,0}^{m-1}$. Since all graphs contributing to  $J_{e_0,0}^{m-1}$ appear exactly once this way, we have

\begin{equation}\label{eq:sum00}
J_{e_0,0}^m = \sum_{1\leq i<j\leq m} (\nu_i+\nu_j+k)\cdot J_{e_0,0}^{m-1} =\left( \sum_{1\leq i<j\leq m} (\nu_i+\nu_j)+\sum_{1\leq i<j\leq m} k\right)\cdot J_{e_0,0}^{m-1}. 
\end{equation}

We observe that \eqref{eq:sum00} is unambiguous (and in fact the last term is well-defined) because by inductive hypothesis the invariants $J_{e_0,0}^{m-1}$ are independent of the negative ramification profiles. 
The first summand in the last term of \eqref{eq:sum00} equals $(m-1)\cdot(\nu_1+\ldots+\nu_m)$, which by \eqref{eq:leakyglobbal} then equals
$(m-1)\cdot(d-k\cdot(m-1))$.
The second summand equals ${{m}\choose{2}}\cdot k$. Thus \eqref{eq:sum00} becomes:
\begin{equation}
 J_{e_0,0}^m = \left((m-1)\cdot(d-k\cdot(m-1))+{{m}\choose{2}}\cdot k\right)\cdot J_{e_0,0}^{m-1},  
\end{equation}
 which is easily seen to agree with \eqref{eq:simplifiedrec} (since $e_0<m-2$ the Kronecker delta equals zero). Formula \eqref{eq:simplifiedrec} immediately establishes the independence of $J^m_{e_0, 0}$ from $\bfnu$.

\end{proof}

\begin{lemma} \label{Lem:genus_0_recursion}
Given $m \geq 3$, $e_0\geq 0$, $e_1>0$ and $e_0+e_1\leq m-2$, we have
\begin{align} \label{Eq:recgzero}
   & J_{e_0,e_1}^m=\\ \nonumber &
\Big(\binom{m-1}{e_1+1}\cdot (k \cdot(e_1+1)+\nu_1)+ \binom{m-2}{e_1}\cdot (d-\nu_1-k\cdot (m-1)) \Big)\cdot J_{e_0,0}^{m-e_1-1}\\& \nonumber
  +\Big(\binom{m-1}{2}\cdot k+(m-2)\cdot (d-\nu_1-k\cdot (m-1))\Big)\cdot J_{e_0,e_1}^{m-1} +\binom{m-2}{e_0}\cdot\delta_{m-2,e_0+e_1}\,,
\end{align}
with the base cases of the recursion again encoded by the Kronecker delta part of the formula. In particular the invariants $J_{e_0, e_1}^m$ are independent of $\nu_2, \ldots, \nu_m$, the ramification conditions at points that do not support a positive power of a $\psi$-class.
\end{lemma}
\begin{proof}
The proof is very similar to the proof of Lemma \ref{Lem:genus_0_recursion0}, so we move along a bit faster. In particular, we omit weaving the inductive proof of independence from $\nu_2, \ldots, \nu_m$ with the establishing of the recursion, as it goes exactly as before and would only burden the interesting part. The base cases are evaluation of two-term monomials in $\psi$-classes on moduli spaces of rational stable curves, which are well-known to produce the above binomial coefficients.

Letting $e_0+e_1<m-2$, and applying the same vertex cutting algorithm as in the proof of Lemma \ref{Lem:genus_0_recursion0}, we must now consider two cases (see Figure \ref{fig:recursionPsi0Psi1}):

\begin{description}
    \item[Case 1] the last vertex is adjacent to the end marked $1$ of weight $\nu_1$, in which case it has to be $e_1+3$-valent. As it is adjacent to a single bounded edge, it must be adjacent to $e_1+1$ more ends. By Remark \ref{rem-localvertexmultgenus0}, its local vertex multiplicity is one, independently of the weight of the other adjacent ends. Let $I\subset \{2\ldots,m\}$ be the subset of additional ends adjacent to the last vertex. We use the notation
    $$\nu_I:= \sum_{i\in I} \nu_i.$$
    The bounded edge which we cut is then of weight
    $$ k \cdot(e_1+1)+\nu_1+\nu_I.$$
The part that remains (to the left of the cut) is a tropical leaky cover satisfying a $\psi_0^{e_0}$ condition on the left end of weight $d$, and has $e_1+1$ fewer right ends. It thus contributes to $J_{e_0,0}^{m-e_1-1}$. Vice versa, one can prolong any such cover by re-attaching the last vertex, the multiplicity changes by a factor of $k \cdot(e_1+1)+\nu_1+\nu_I$. Altogether, this first case contributes
\begin{equation}
\label{eq:sumoversubsets}
\sum_{\tiny{
\begin{array}{c}
  I\subset \{2,\ldots,m\}     \\
  | I| = e_1+1    
\end{array}
}}
\big(k \cdot(e_1+1)+\nu_1+\nu_I\big) \cdot J_{e_0,0}^{m-e_1-1}. 
\end{equation}

To simplify \eqref{eq:sumoversubsets}, we rewrite it as follows:
\begin{equation*}
\sum_{\tiny{
\begin{array}{c}
  I\subset \{2,\ldots,m\}     \\
  | I| = e_1+1    
\end{array}
}}
\big(k \cdot(e_1+1)+\nu_1\big) \cdot J_{e_0,0}^{m-e_1-1}+
\sum_{\tiny{
\begin{array}{c}
  I\subset \{2,\ldots,m\}     \\
  | I| = e_1+1    
\end{array}
}}
\nu_I \cdot J_{e_0,0}^{m-e_1-1}
. 
\end{equation*}
The  first summation consists of $\binom{m-1}{e_1+1}$  terms which do not depend on $I$, and therefore add up to 
\begin{equation}\label{eq:firstterm}
\binom{m-1}{e_1+1}
\cdot \big(k \cdot(e_1+1)+\nu_1\big)\cdot   J_{e_0,0}^{m-e_1-1}   
\end{equation}

For the second summation, it is a simple exercise in symmetric functions that 
\begin{equation}\label{Eq:symmsumm}
  \sum_{\tiny{
\begin{array}{c}
  I\subset \{2,\ldots,m\}     \\
  | I| = e_1+1    
\end{array}
}}
\nu_I  = \binom{m-2}{e_1} \cdot\sum_{i=2}^{m}\nu_i.
\end{equation}
From the $k$-leaky balancing condition \eqref{eq:leakyglobbal} we have:
\begin{equation}\label{eq:leakybalglobal}
\sum_{i=2}^m \nu_i   =
(d-\nu_1-k\cdot (m-1))
\end{equation}
Combining \eqref{Eq:symmsumm}
and \eqref{eq:leakybalglobal} one obtains that the second summation equals:
\begin{equation}
\label{eq:secondumm}
 \binom{m-2}{e_1} 
\cdot(d-\nu_1-k\cdot (m-1))\cdot  J_{e_0,0}^{m-e_1-1}. 
\end{equation}
Adding \eqref{eq:firstterm}, \eqref{eq:secondumm} it is immediate to see that \eqref{eq:sumoversubsets} agrees with the first line of the recursion in \eqref{Eq:recgzero}. 


\item[Case 2] If the last vertex is not adjacent to the end of weight $\nu_1$, it must be $3$-valent and have multiplicity $1$. If it is adjacent to the ends of weight $\nu_i$ and $\nu_j$, where $i,j\neq 1$, the weight of its adjacent bounded edge is $\nu_i+\nu_j+k$. Cutting this edge, the remainder of the graph (to the left of the cut) contributes to $J_{e_0,e_1}^{m-1}$.
Similarly to before, to any  contributing graph one may attach a tripod to undo the cut. Thus, the second case contributes

\begin{equation}\label{eq:sumcase2}
\sum_{\{i,j\}\subseteq\{2,\ldots,m\}} (\nu_i+\nu_j+k)\cdot J_{e_0,e_1}^{m-1}. 
\end{equation} 

Expression \eqref{eq:sumcase2} is linear symmetric in $\nu_2, \ldots, \nu_m$. Every $\nu_i$ appears in $m-2$ terms, thus we obtain 
\begin{align} \label{eq:Secontermofrec}\nonumber&\Big(\binom{m-1}{2}\cdot k+(m-2)\cdot (\sum_{i=2}^m \nu_i)\Big)\cdot J_{e_0,e_1}^{m-1}=\\&
\Big(\binom{m-1}{2}\cdot k+(m-2)\cdot (d-\nu_1-k\cdot (m-1))\Big)\cdot J_{e_0,e_1}^{m-1}.
\end{align}

\end{description}

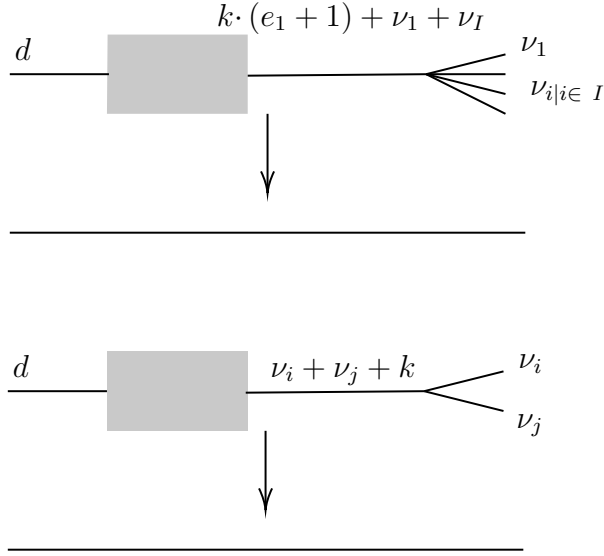
\begin{figure}
    \centering

\tikzset{every picture/.style={line width=0.75pt}} 

\begin{tikzpicture}[x=0.75pt,y=0.75pt,yscale=-1,xscale=1]

\draw  [draw opacity=0][fill={rgb, 255:red, 155; green, 155; blue, 155 }  ,fill opacity=0.54 ] (199.33,320) -- (270,320) -- (270,360) -- (199,360) -- cycle ;
\draw    (270,340.75) -- (360,340) ;
\draw    (360,340) -- (400,330) ;
\draw    (360,340) -- (400,350) ;
\draw    (150,340) -- (200,340) ;
\draw    (280,360) -- (280,398) ;
\draw [shift={(280,400)}, rotate = 270] [color={rgb, 255:red, 0; green, 0; blue, 0 }  ][line width=0.75]    (10.93,-3.29) .. controls (6.95,-1.4) and (3.31,-0.3) .. (0,0) .. controls (3.31,0.3) and (6.95,1.4) .. (10.93,3.29)   ;
\draw    (150,420) -- (410,420) ;
\draw  [draw opacity=0][fill={rgb, 255:red, 155; green, 155; blue, 155 }  ,fill opacity=0.54 ] (199,480) -- (270,480) -- (270,520) -- (199,520) -- cycle ;
\draw    (269,500.75) -- (359,500) ;
\draw    (359,500) -- (399,490) ;
\draw    (359,500) -- (399,510) ;
\draw    (149,500) -- (199,500) ;
\draw    (279,520) -- (279,558) ;
\draw [shift={(279,560)}, rotate = 270] [color={rgb, 255:red, 0; green, 0; blue, 0 }  ][line width=0.75]    (10.93,-3.29) .. controls (6.95,-1.4) and (3.31,-0.3) .. (0,0) .. controls (3.31,0.3) and (6.95,1.4) .. (10.93,3.29)   ;
\draw    (149,580) -- (409,580) ;
\draw    (360,340) -- (400,340) ;
\draw    (360,340) -- (400,360) ;

\draw (151,320) node [anchor=north west][inner sep=0.75pt]   [align=left] {$\displaystyle d$};
\draw (406,320) node [anchor=north west][inner sep=0.75pt]   [align=left] {$\displaystyle \nu_{1}$};
\draw (411,342) node [anchor=north west][inner sep=0.75pt]   [align=left] {$\displaystyle \nu_{i|i\in \ I}$};
\draw (281,320) node [anchor=north west][inner sep=0.75pt]   [align=left] {$ $};
\draw (150,480) node [anchor=north west][inner sep=0.75pt]   [align=left] {$\displaystyle d$};
\draw (405,480) node [anchor=north west][inner sep=0.75pt]   [align=left] {$\displaystyle \nu_{i}$};
\draw (404,510) node [anchor=north west][inner sep=0.75pt]   [align=left] {$\displaystyle \nu_{j}$};
\draw (280,480) node [anchor=north west][inner sep=0.75pt]   [align=left] {$\displaystyle \nu_{i} +\nu_{j} +k$};
\draw (253,302) node [anchor=north west][inner sep=0.75pt]   [align=left] {$\displaystyle k\cdotp ( e_{1} +1) +\nu_{1} +\nu_{I}$};

\end{tikzpicture}

    \caption{Sketch for the recursion for the invariants $J_{e_0,e_1}(m)$. The tropical leaky covers  continue in some way within the grey boxes, we only pay attention to the behavior of the last vertex.}
    \label{fig:recursionPsi0Psi1}
\end{figure}

Adding the contributions \eqref{eq:firstterm}, \eqref{eq:secondumm} and \eqref{eq:Secontermofrec} from the two cases, we obtain recursion \eqref{Eq:recgzero}.


\end{proof}

The recursions from Lemmas \ref{Lem:genus_0_recursion0}, \ref{Lem:genus_0_recursion}  are used to prove the formulas from Theorem \ref{Thm:genus0} .

\begin{proof}[Proof of Theorem \ref{Thm:genus0} \eqref{eqn:J_e0_0}, case $e_1=0$: ]
For $e_1=0$, we use the recursion from Lemma \ref{Lem:genus_0_recursion0}. 
Having already established the base cases, we show that for $e_0<m-2$ the right hand side of equation \eqref{eqn:J_e0_0} satisfies the recursion from Lemma \ref{Lem:genus_0_recursion}, then the statement follows. 

Inserting
$J_{e_0,0}^{m-1}= \frac{(m-2)!}{(e_0+1)!} \prod_{j=e_0+1}^{m-3}(d-j\frac{k}{2})$ into the right hand side of \eqref{eq:simplifiedrec}, we obtain

\begin{align*}
&(m-1)\cdot \Big(d-(m-2)\frac{k}{2}\Big)
\frac{(m-2)!}{(e_0+1)!} \cdot \Big( d-(e_0+1)\frac{k}{2}\Big)\cdot \ldots \cdot \Big( d-(m-3)\frac{k}{2}\Big)=\\& \frac{(m-1)!}{(e_0+1)!} \cdot \Big( d-(e_0+1)\frac{k}{2}\Big)\cdot \ldots \cdot \Big( d-(m-2)\frac{k}{2}\Big)=J_{e_0,0}^m\end{align*}
as required.

\end{proof}
\begin{proof}[Proof of Theorem \ref{Thm:genus0} \eqref{eqn:J_0_e1}, case $e_0=0$]
The base case is $J^m_{0,e_1} = 1$ when $e_1 = m-2$, which agrees with the right hand side of $\eqref{eqn:J_0_e1}$ evaluated at $e_1 = m-2$.

To show that the right hand side of \eqref{eqn:J_0_e1} satisfies recursion \eqref{Eq:recgzero}, we introduce temporarily the notation:
\begin{equation}
    \widehat{J}^m_{0,e_1} := {J}^m_{0,e_1}\frac{(e_1+1)!}{(m-2)!} = F^m+G^m+H^m,
\end{equation}
where:\\
\resizebox{\textwidth}{!}{$
F^m:= e_1 \cdot \frac{\prod_{j=e_1+1}^{m-1} (d-\nu_1-j\cdot\frac{k}{2}) }{-\nu_1-e_1\cdot\frac{k}{2}}, 
 G^m:= e_1 \cdot \frac{ - \prod_{j=1}^{m-e_1-1} (d-j\cdot\frac{k}{2})}{-\nu_1-e_1\cdot\frac{k}{2}},  
  H^m:=  (m-e_1-1) \cdot \prod_{j=1}^{m-e_1-2} \left(d-j\cdot\frac{k}{2}\right).$}

Plugging $J_{0,0}^{m-e_1-1} = (m-e_1-2)! \prod_{j = 1}^{m-e_1-2} \left(d-j\frac{k}{2}\right) $ and doing some elementary factorizations, one sees that \eqref{Eq:recgzero} is equivalent to:
\begin{equation}\label{eq:chefatica}
 \widehat{J}^m_{0,e_1} = ((d-\nu_1)(e_1+1)+\nu_1(m-1))\cdot \prod_{j = 1}^{m-e_1-2}\left(d-j\cdot\frac{k}{2}\right)  +
 (d-\nu_1-(m-1)\cdot\frac{k}{2})\cdot \widehat{J}^{m-1}_{0,e_1}.
\end{equation}
Consider the second summand of
\eqref{eq:chefatica}. It is immediate that:
\begin{equation}\label{eq:cf1}
    F^m = \left(d-\nu_1-(m-1)\cdot\frac{k}{2}\right)\cdot F^{m-1}.  
\end{equation}
By adding and subtracting  $e_1\cdot \frac{k}{2}$ one sees:
\begin{equation}\label{eq:cf2}
  G^m =  \left(d-\nu_1-(m-1)\cdot\frac{k}{2}\right)\cdot G^{m-1} +\left(\nu_1+ e_1\cdot\frac{k}{2}\right)\cdot G^{m-1} .  
\end{equation}
One can conclude the proof by showing that
\begin{align}\label{eq:cf3}
H^m &=  ((d-\nu_1)(e_1+1)+\nu_1(m-1))\cdot \prod_{j = 1}^{m-e_1-2}\left(d-j\cdot\frac{k}{2}\right)  +
 (d-\nu_1-(m-1)\cdot\frac{k}{2})\cdot H^{m-1} \nonumber\\
 & -\left(\nu_1+ e_1\cdot\frac{k}{2}\right)\cdot G^{m-1}, 
\end{align}
since \eqref{eq:chefatica} is obtained by adding  \eqref{eq:cf1}, \eqref{eq:cf2}, \eqref{eq:cf3}.
Factoring out $\prod_{j = 1}^{m-e_1-2}\left(d-j\cdot\frac{k}{2}\right)$ from all terms, \eqref{eq:cf3} reduces to a quadratic expression which is easily verified. 
\end{proof}

\subsection{Polynomiality}

We turn our attention to general genus zero invariants. While we can't produce explicit general formulas, we show that for fixed degree of the descendant insertions, the invariants are polynomial in the ramification profiles and the leaking.  
\begin{lemma} \label{lem:reggenzerogenerale}
 Fix  $\mathbf{e}=(e_0, \ldots, e_m) \in \mathbb{Z}^n$ with $e_0\geq0$, $e_1>0,\ldots,e_u>0$ and $e_{u+1}=\ldots=e_{m}=0$, let $c=m-2-e_0-\ldots-e_u>0$.   We have
\begin{align}\label{eq:onemorerecursion}
\mathsf{H}_0((d, - \bfnu), \mathbf{e}) = \sum_{(I,J)} F_{I,J} \cdot \binom{e_I}{e_{i_1},\ldots,e_{i_{|I|}}}\cdot  (\nu_I+\nu_J+(e_I+1)\cdot k)
      \end{align}
  where   $I\subset\{1,\ldots,u\}$,  $J\subset \{u+1\ldots,m\}$ such that $|I|+|J| = \sum_{i\in I} e_i+2$, $e_I = \sum_{i\in I} e_i$ (similarly for $\nu_I,\nu_J$) and
  \begin{align}\label{al:cutpiece}
    F_{I,J} = \int_{\logDR_0(\bfnu\setminus \{\nu_l\;|\;l\in I\cup J\}) \cup \{\nu_I+\nu_J+(e_I+1)\cdot k\})} \prod_{j\in \{1,\ldots,u\}\setminus I} \psi_j^{e_j}\cdot \branch_{c-1}. 
\end{align}
\end{lemma}
\begin{proof}
    We obtain this recursion by cutting all tropical leaky covers which contribute to $\mathsf{H}_0((d, - \bfnu), \mathbf{e})$ right before their last vertex. The right pointing ends attached to such a vertex identify a subset $I$ of $\{1, \ldots, u\}$ and a subset $J$ of $\{u+1, \ldots, m\}$ such that $|I|+|J| = \sum_{i\in I} e_i+2$, and vice-versa for any pair of subsets we can find graphs with a vertex of that type, which we call $V_{I,J}$. The unique left pointing end of $V_{I,J}$ has weight $(\nu_I+\nu_J+(e_I+1)\cdot k)$, where $\nu_I = \sum_{i\in I} \nu_i$ and similarly for $\nu_J$.
    Finally, consider all the graphs that contain a vertex of type $V_{I,J}$: the remaining parts of the graphs to the left of the cut are precisely the leaky tropical covers contributing to \eqref{al:cutpiece}, thus concluding the proof of the lemma.
\end{proof}

\begin{proposition} \label{Prop:onepart_genus0_structure}

With notation as above,  the invariant
\[
\mathsf{H}_0((d, - \bfnu), \mathbf{e}) = \int_{\logDR_0(d, -\bfnu)} \psi_0^{e_0}\psi_1^{e_1}\ldots\psi_u^{e_u}\cdot \branch_{c}\,
\]
is  a polynomial of degree $c$ in $d$, $k$, and $\nu_1,\ldots,\nu_u$. In particular, it is independent of the variables $\nu_{u+1}, \ldots, \nu_{m}$ associated to markings without $\psi$-insertions.
\end{proposition}

\begin{proof}

We prove this by induction on $c=m-2-e_0-\ldots-e_u$. If $c=0$, the integral above reduces to the known integral
\[
\int_{\overline{\mathcal{M}}_{0,m+1}} \psi_0^{e_0} \cdots \psi_u^{e_u} =  \binom{m-2}{e_0, e_1, \ldots, e_u}\,,
\]
which is constant in $d,k$ and the entries of $\mathbf{\nu}$. On the tropical side, this corresponds to a single $k$-leaky cover with just one vertex and thus contributing its multiplicity. 

Now consider the case $c>0$. By Lemma \ref{lem:reggenzerogenerale}, $ \mathsf{H}_0(d, - \bfnu)$ may be written as a finite sum, where the indexing set $(I,J)$ of the summation does not depend on $d, k$ or any $\nu_i$. Each summand is the product of three terms:  the first, $F_{I,J}$ is an invariant with branch degree equal to $c-1$, therefore by the inductive hypothesis it is a polynomial of degree $c-1$ in (a subset of) the variables $d, k$ and $\nu_1, \ldots, \nu_u$: note that the end of weight  $(\nu_I+\nu_J+(e_I+1)k)$ does not support any descendant insertion, and therefore the invariant does not depend on this weight.
The second factor is a multinomial coefficient depending only on the descendant vector $\mathbf{e}$, therefore it is constant in our variable of interest. The last factor is linear in $k$ and all the $\nu_i$'s, so we can already see that $ \mathsf{H}_0(d, - \bfnu)$ is a polynomial of degree $c$ in $d, k$ and the end weights. We must now show that we can eliminate the dependence on the $\nu_j$ with $j\in J$.
We observe that by the inductive hypothesis $F_{I,J}$ only depends on the subset $I$ and not on $J$, so we write $F_I = F_{I,J}$ for any valid choice of $J$. We can hence rewrite \eqref{eq:onemorerecursion} as 
\begin{equation} \label{eqn:H_genus0_more_general}
    \mathsf{H}_0((d, - \bfnu), \mathbf{e}) = \sum_{I} F_I \cdot \binom{e_I}{e_{i_1},\ldots,e_{i_{|I|}}}\cdot \sum_{J}  \cdot (\nu_I+\nu_J+(e_I+1)\cdot k)\,,
\end{equation}
and we focus on the second summation.
The summand $\nu_I+(e_I+1)\cdot k$ is the same for each choice of $J$, so in total we obtain it 
$ \binom{m-u}{e_I+2-|I|}$
times. For every $u+1\leq j\ \leq m$,  $\nu_j$ appears in the summation
$\binom{m-u-1}{e_I+1-|I|}$
 times, so we obtain the sum 
 \[\nu_{u+1}+\ldots+\nu_m = d-(m-1)k-\nu_1-\ldots-\nu_u\]
 multiplied with this binomial coefficient.
Substituting in \eqref{eqn:H_genus0_more_general}, we obtain
\begin{align*}
 \mathsf{H}_0((d, - \bfnu), \mathbf{e}) = \sum_{I} F_I \cdot \binom{e_I}{e_{i_1},\ldots,e_{i_{|I|}}}&\cdot \Bigg(\binom{m-u}{e_I+2-|I|}
\cdot (\nu_I+(e_I+1)k) \;\;+\\&
\binom{m-u-1}{e_I+1-|I|} \cdot (d-(m-1)k-\nu_1-\ldots-\nu_u)\Bigg),
\end{align*}
which is a polynomial of degree $c$ in $d, k, \nu_1, \ldots ,\nu_u$.
\end{proof}

\section{Higher genus}

In this section we turn our attention to higher genus invariants. In this case we still obtain a recursion  by excising the last vertex of tropical covers, but the presence of positive genus allows for more possible terms to the recursion. We  then  provide some formulas in genus one, and some in arbitrary genus for $k=0$.

\subsection{Recursion for one-part leaky double Hurwitz numbers}

We formulate a recursion for leaky covers which allows for arbitrary genus, but restricted to the case $\mathbf{e}=\mathbf{0}.$
\begin{lemma}\label{prop-rechighergenus}
For genus $g \geq 0$ and $m \geq 1$, the one-part $k$-leaky double Hurwitz numbers $\mathsf{H}_g(d, \bfnu)$ 
satisfy the following recursion: 
\begin{align}
\mathsf{H}_g(d, \bfnu) &= \sum_{1\leq i<j\leq m} (\nu_i+\nu_j+k) \cdot \mathsf{H}_g(d, \bfnu \setminus \{\nu_i,\nu_j\} \cup \{\nu_i+\nu_j+k\}) \\
&+ \frac{1}{2} \sum_{a+b= \nu_i+k} a \cdot b \cdot \mathsf{H}_{g-1}(d, \bfnu \setminus \{\nu_i\} \cup \{a, b\}) \\
&- \frac{k}{24} \cdot \mathsf{H}_{g-1}(d, \bfnu \cup \{k\}) \\
&+ \frac{1}{6} \sum_{a+b+c=k} a \cdot b \cdot c \cdot \mathsf{H}_{g-2}(d, \bfnu \cup \{a,b,c\}),
\end{align}
where $a,b,c$ are intended to be non-negative integers; we have as base cases:
\begin{align}\label{eq:basecase}
\mathsf{H}_0(d, \bfnu) &= (m-1)! \cdot \prod_{j=1}^{m-2} \left(d - j \cdot \frac{k}{2}\right)
\end{align}
\end{lemma}

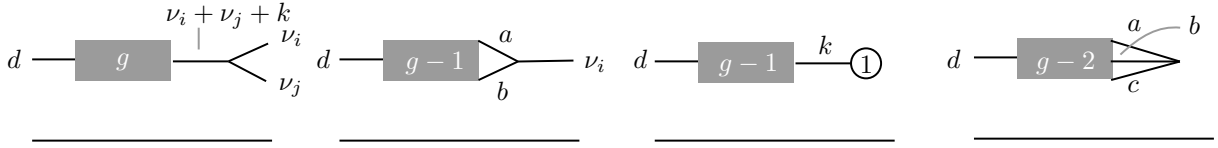
\begin{figure}
    \centering

\tikzset{every picture/.style={line width=0.75pt}} 

\begin{tikzpicture}[x=0.75pt,y=0.75pt,yscale=-1,xscale=1]

\draw    (21.14,250.43) -- (142.14,250.43) ;
\draw    (21.14,209.43) -- (43.14,209.43) ;
\draw  [draw opacity=0][fill={rgb, 255:red, 155; green, 155; blue, 155 }  ,fill opacity=1 ] (43,200) -- (91.14,200) -- (91.14,220.43) -- (43,220.43) -- cycle ;
\draw    (92.14,210.43) -- (120.14,210.43) ;
\draw    (120.14,210.43) -- (140.14,201.43) ;
\draw    (139.14,220.43) -- (120.14,210.43) ;
\draw    (176.14,250.43) -- (297.14,250.43) ;
\draw    (176.14,209.43) -- (198.14,209.43) ;
\draw  [draw opacity=0][fill={rgb, 255:red, 155; green, 155; blue, 155 }  ,fill opacity=1 ] (198,200) -- (246.14,200) -- (246.14,220.43) -- (198,220.43) -- cycle ;
\draw    (294.15,210.02) -- (266.15,210.5) ;
\draw    (266.15,210.5) -- (246.31,219.84) ;
\draw    (246.14,200) -- (266.15,210.5) ;
\draw    (335.14,250.43) -- (456.14,250.43) ;
\draw    (335.14,210.43) -- (357.14,210.43) ;
\draw  [draw opacity=0][fill={rgb, 255:red, 155; green, 155; blue, 155 }  ,fill opacity=1 ] (357,201) -- (405.14,201) -- (405.14,221.43) -- (357,221.43) -- cycle ;
\draw    (406.14,211.43) -- (434.14,211.43) ;
\draw    (496.14,249.43) -- (617.14,249.43) ;
\draw    (496.14,208.43) -- (518.14,208.43) ;
\draw  [draw opacity=0][fill={rgb, 255:red, 155; green, 155; blue, 155 }  ,fill opacity=1 ] (518,199) -- (566.14,199) -- (566.14,219.43) -- (518,219.43) -- cycle ;
\draw   (434.14,211.43) .. controls (434.14,207.25) and (437.53,203.86) .. (441.71,203.86) .. controls (445.9,203.86) and (449.29,207.25) .. (449.29,211.43) .. controls (449.29,215.61) and (445.9,219) .. (441.71,219) .. controls (437.53,219) and (434.14,215.61) .. (434.14,211.43) -- cycle ;
\draw    (600,210.5) -- (565.31,219.84) ;
\draw    (565.14,200) -- (600,210.5) ;
\draw    (600,210.5) -- (565.14,210.43) ;
\draw [color={rgb, 255:red, 155; green, 155; blue, 155 }  ,draw opacity=1 ][line width=0.75]    (105.14,193.43) -- (105.14,204.43) ;
\draw [color={rgb, 255:red, 155; green, 155; blue, 155 }  ,draw opacity=1 ]   (570,207) .. controls (580.14,197.43) and (589.14,192.43) .. (600.14,193.43) ;

\draw (437,206) node [anchor=north west][inner sep=0.75pt]  [font=\footnotesize]  {$1$};
\draw (62,205) node [anchor=north west][inner sep=0.75pt]  [font=\footnotesize,color={rgb, 255:red, 255; green, 255; blue, 255 }  ,opacity=1 ]  {$g$};
\draw (208,204) node [anchor=north west][inner sep=0.75pt]  [font=\footnotesize,color={rgb, 255:red, 255; green, 255; blue, 255 }  ,opacity=1 ]  {$g-1$};
\draw (366,205) node [anchor=north west][inner sep=0.75pt]  [font=\footnotesize,color={rgb, 255:red, 255; green, 255; blue, 255 }  ,opacity=1 ]  {$g-1$};
\draw (526,204) node [anchor=north west][inner sep=0.75pt]  [font=\footnotesize,color={rgb, 255:red, 255; green, 255; blue, 255 }  ,opacity=1 ]  {$g-2$};
\draw (87,179.4) node [anchor=north west][inner sep=0.75pt]  [font=\footnotesize]  {$\nu _{i} +\nu _{j} +k$};
\draw (145,193.4) node [anchor=north west][inner sep=0.75pt]  [font=\footnotesize]  {$\nu _{i}$};
\draw (144,217.4) node [anchor=north west][inner sep=0.75pt]  [font=\footnotesize]  {$\nu _{j}$};
\draw (7,203) node [anchor=north west][inner sep=0.75pt]  [font=\footnotesize]  {$d$};
\draw (163,203) node [anchor=north west][inner sep=0.75pt]  [font=\footnotesize]  {$d$};
\draw (323,203) node [anchor=north west][inner sep=0.75pt]  [font=\footnotesize]  {$d$};
\draw (481,203) node [anchor=north west][inner sep=0.75pt]  [font=\footnotesize]  {$d$};
\draw (298,206.4) node [anchor=north west][inner sep=0.75pt]  [font=\footnotesize]  {$\nu _{i}$};
\draw (416,196) node [anchor=north west][inner sep=0.75pt]  [font=\footnotesize]  {$k$};

\draw (255,193.4) node [anchor=north west][inner sep=0.75pt]  [font=\footnotesize]  {$a$};
\draw (572,188) node [anchor=north west][inner sep=0.75pt]  [font=\footnotesize]  {$a$};
\draw (254,219.4) node [anchor=north west][inner sep=0.75pt]  [font=\footnotesize]  {$b$};
\draw (603,186) node [anchor=north west][inner sep=0.75pt]  [font=\footnotesize]  {$b$};
\draw (572,218.57) node [anchor=north west][inner sep=0.75pt]  [font=\footnotesize]  {$c$};

\end{tikzpicture}

    \caption{The four possible local structures of leaky covers at the last vertex.}
    \label{fig:lostru}
\end{figure}
\begin{proof}
The proof is analogous to the recursion in Lemma \ref{Lem:genus_0_recursion}.  We cut the last vertex of a $k$-leaky tropical cover, compute the contribution of this cut part, and notice that the remaining part of the cover contributes to an  invariant where the quantity $(2g-1+m)$ has decreased by one. The four terms in this recursion correspond to different local structures at the last vertex of a tropical $k$-leaky cover, as illustrated in Figure \ref{fig:lostru}:
\begin{enumerate}
\item The last vertex carries two markings $\nu_i$ and $\nu_j$, connected to the rest of the curve by one bounded edge.
\item The last vertex carries one marking $\nu_i$, connected to the rest of the curve by two bounded edges, which lowers the genus of the remaining part by  $1$.
\item The last vertex forms a cul-de-sac (a one-valent genus one vertex with one bounded edge and no markings), which lowers the genus by $1$.
\item The last vertex has three bounded edges connecting to the rest of the curve, which lowers the genus by $2$.
\end{enumerate}

In this formula, $k$ is determined by the balancing condition $d-\sum_{i=1}^m \nu_i = k \cdot(2g+m-1)$.
In the first case, the last vertex has multiplicity one and the bounded edge which we cut has weight $\nu_i+\nu_j+k$. The remaining cover after cutting this bounded edge still has genus $g$. Instead of the two right ends of weights $\nu_i,\nu_j$ it has one right end of weight $\nu_i+\nu_j+k$.

In the second case, the last vertex is still $3$-valent, but we cut two bounded edges. If we set the weight of one of these edges to be $a$, then by the leaky balancing condition, the second has weight $b = \nu_i+k-a$. Summing over all $a,b$ such that $a+b = \nu_i+k$ we double count isomorphic leaky covers that appear. For that reason, we multiply with a factor of $\frac{1}{2}$. 

In the third case, the vertex we cut off has a local multiplicity of $-\frac{1}{24}$ (see e.g.\ Example 6.1 \cite{CMS24}). The genus of the remaining cover is $g-1$, and it has a new end, of weight $k$.

For the fourth case, the local vertex multiplicity is one. We cut three edges whose weights we can denote by $a$, $b$, and, by leaky balancing, $c = k-a-b$. The factor of $\frac{1}{6}$ is again obtained by the $S_3$ symmetry in the roles of the three edges, which leads to overcounting covers by a factor of $6$ when summing over all $a,b,c$ such that $a+b+c = k$. 
\end{proof}

The recursion from Proposition \ref{prop-rechighergenus} allows us to generalize the result from \cite[Corollary 3.3]{GJV05} on genus one one-part double Hurwitz numbers to arbitrary $k$ and prove Theorem \ref{thm:genusone}.


\begin{proof}[Proof of Theorem \ref{thm:genusone}]
We apply the recursion from Proposition \ref{prop-rechighergenus} in the case  $g=1$, inserting the formula from \eqref{eqn:GJV_genus_1} and the base case \eqref{eq:basecase}; the right hand side of the recursion is:

\begin{align*}\hspace{-1cm} \label{eq:g1recapply}
   & \sum_{\{i,i'\}\subseteq \{1,\ldots,m\}} (\nu_i+\nu_{i'}+k)\cdot \frac{m!}{24}\cdot \prod_{j=2}^{m-1} \left(d - j\cdot \frac{k}{2}\right)
    \cdot\left[(d-k)^3 - \left(d-\frac{k}{2}\right)\right. \\ & \left.
    \cdot\left(1+   2\cdot \sum_{s\not\in \{i,i'\}}(\nu_i+\nu_{i'}+2k)\cdot (\nu_s+k)
+2\cdot \sum_{\tiny{\begin{array}{c}
       \{s,t\}\cap\{i,i'\} = \phi   
    \end{array}}}(\nu_s+k)\cdot(\nu_t+k)
\right)\right]
\\& +
\frac{1}{12}\cdot \left(\sum_{i= 1}^m ((\nu_i+k)^3-(\nu_i+k))-\frac{k}{2}\right)\cdot m!\cdot \prod_{j=1}^{m-1} \left(d - j\cdot \frac{k}{2}\right)\,,
\end{align*}
where we used that
\[
\frac{1}{2} \sum_{a + b = \nu_i + k} a \cdot b = \frac{1}{12}((\nu_i+k)^3 - (\nu_i+k))\,.
\]
We factor $\frac{m!}{24}\cdot \prod_{j=2}^{m-1} \left(d - j\cdot \frac{k}{2}\right)$ from 
the recursion above, and make the substitution $x_i = \nu_i+k$, to obtain
\begin{align} 
\nonumber
& \sum_{\{i,i'\}\subseteq \{1,\ldots,m\}}(x_i+x_{i'}-k)  \left[(d-k)^3-\left(d-\frac{k}{2}\right) \right.\\ & \left.\cdot\left(1+2\cdot\sum_{s\not\in \{i,i'\}} (x_i+x_{i'})x_s+ 2\cdot\sum_{\tiny{\begin{array}{c}
       \{s,t\}\cap\{i,i'\} = \phi   
    \end{array}}} x_sx_t\right)\right] \label{eq:recwrittenwithxi}
    \\& + 2\cdot \left((d-k)^3 -3(d-k)\cdot \sigma_2+3\sigma_3-(d-k)-\frac{k}{2}\right) \cdot\left(d - \frac{k}{2}\right), \nonumber
\end{align}
where we have manipulated the third line by  expressing the sum of cubes in the elementary symmetric functions in the $x_i$'s
$$
\sum_{i=1}^m x_i^3 = \sigma_1^3 - 3\sigma_1\sigma_2+3\sigma_3
$$
and using the relation $\sigma_1 = \sum_{i=1}^m x_i = d-k$.

We must prove that \eqref{eq:recwrittenwithxi} 
equals formula \eqref{eqn:GJV_genus_1} after factoring out $\frac{m!}{24}\cdot \prod_{j=2}^{m-1} \left(d - j\cdot \frac{k}{2}\right)$; adopting the symmetric function notation this is:
\begin{equation}\label{eq:whatmustbecomparredto}
    (m+1)\cdot\left(d - m\cdot\frac{k}{2}\right) \cdot \left[(d-k)^3 - \left(d -  \frac{k}{2}\right)\cdot\left(1+2\sigma_2
    \right)
    \right].
\end{equation}
Both \eqref{eq:recwrittenwithxi} and \eqref{eq:whatmustbecomparredto} are polynomials in $d, k$ and the $x_i$'s with a homogeneous term of degree $4$ and one of degree $2$. We begin by comparing the degree $2$ terms, and observe that they share a common factor of  $\left(d - \frac{k}{2}\right)$. The remaining linear part in \eqref{eq:recwrittenwithxi} is 
\begin{align}
 &-\sum_{\{i,i'\}\subseteq \{1,\ldots,m\}}(x_i+x_{i'}-k)  - 2d +k \nonumber \\
 &= - (m-1)\cdot(d-k) +{{m}\choose{2}}\cdot k -2d+k   \nonumber \\
 &= -(m+1)\cdot \left(d - m\cdot \frac{k}{2}\right),
\end{align}
also agreeing with the remaining linear part in \eqref{eq:whatmustbecomparredto}.

Next we observe there are parts of the degree $4$ terms which are divisible by $(d-k)^3$. Factoring this term out, in \eqref{eq:recwrittenwithxi} we have:
\begin{align}
 &\sum_{\{i,i'\}\subseteq \{1,\ldots,m\}}(x_i+x_{i'}-k)  +2d -k \nonumber \\
 &= -(m+1)\cdot \left(d - m\cdot \frac{k}{2}\right),
\end{align}
which agrees with the corresponding part in \eqref{eq:whatmustbecomparredto}.

Factoring out  $2\cdot \left(d - \frac{k}{2}\right)$ from all remaining terms, it remains to show that
\begin{align} \label{eq:recwrittenwithxiremaining}
\nonumber
 -\sum_{i\not=i'\in \{1,\ldots,m\}}(x_i+x_{i'}-k)  \cdot\left(\sum_{s\not\in \{i,i'\}} (x_i+x_{i'})x_s+\sum_{\tiny{\begin{array}{c}
       \{s,t\}\cap\{i,i'\} = \phi   
    \end{array}}} x_sx_t\right) \\
     -3\cdot(d-k)\cdot\sigma_2+3\sigma_3  =  - (m+1)\cdot\left(d - m\cdot\frac{k}{2}\right)\cdot \sigma_2.
\end{align}

The first line on \eqref{eq:recwrittenwithxiremaining} is a symmetric polynomial in the $x_i$'s, which we manipulate to express as polynomial in the elementary symmetric functions $\sigma_1, \sigma_2, \sigma_3$.
For the summands that are a multiple of $k$, we have
\begin{equation}\label{eq:yesk}
 k\cdot \sum_{\{i,i'\}\subseteq \{1,\ldots,m\}} \left(\sum_{s\not\in \{i,i'\}} (x_i+x_{i'})x_s+\sum_{\tiny{\begin{array}{c}
       \{s,t\}\cap\{i,i'\} = \phi   
    \end{array}}} x_sx_t\right)  = k\cdot \frac{(m-2)\cdot(m+1)}{2}\cdot  \sigma_2,
\end{equation}
since there are $2\cdot(m-2)$ ways for a monomial of the form $x_ax_b$ to appear when $a$ or $b$ are in $\{i,i'\}$ and ${{m-2}\choose{2}}$ ways when $\{a,b\}\cap\{i,i'\} = \phi$.

For the remaining terms, we have:
\begin{align} \label{eq:nok}
  &-\sum_{\{i,i'\}\subseteq \{1,\ldots,m\}}\cdot\left(\sum_{s\not\in \{i,i'\}} (x_i+x_{i'})^2x_s+\sum_{\tiny{\begin{array}{c}
       \{s,t\}\cap\{i,i'\}    = \phi   
    \end{array}}} (x_i+x_{i'})  x_sx_t\right)  \nonumber\\
     &=- (m-2 )\cdot \sum_{a,b = 1}^m x_a^2x_b - \left(2\cdot 3+ 3\cdot(m-3)
     \right) \cdot \sigma_3 \nonumber \\
     & =- (m-2)\cdot  \sigma_1\sigma_2 + 3\cdot((m-2)-2-(m-3))\cdot \sigma_3 \nonumber \\
     & = -(m-2)\cdot  \sigma_1\sigma_2 - 3\cdot \sigma_3, 
\end{align}
where from the second to the third line we have applied the identity:
$$
 \sum_{a,b = 1}^m x_a^2x_b = \sigma_1\sigma_2 - 3 \sigma_3.
$$
Substituting $\sigma_1 = d-k$ and the results of \eqref{eq:yesk}, \eqref{eq:nok} in  the left hand side of \eqref{eq:recwrittenwithxiremaining} we have:
\begin{align}
 &-(m-2)\cdot  (d-k)\cdot\sigma_2 - 3\cdot \sigma_3 + k\cdot \frac{(m-2)\cdot(m+1)}{2}\cdot  \sigma_2   -3\cdot(d-k)\cdot\sigma_2+3\sigma_3 \nonumber \\ &= -(m+1)\cdot d\cdot \sigma_2 +\frac{m \cdot (m+1)}{2}\cdot k\cdot \sigma_2  \nonumber \\
 &=  - (m+1)\cdot\left(d - m\cdot\frac{k}{2}\right)\cdot \sigma_2,
\end{align}
which agrees with the right hand side of \eqref{eq:recwrittenwithxiremaining}, thus concluding the proof.
\end{proof}






\subsection{Proof of Theorem \ref{Thm:psi0_e_general_g}}
In this section we prove Theorem \ref{Thm:psi0_e_general_g} by a different type of recursive structure, obtained by exchanging a $\psi$-class with a boundary divisor. We specialize to the case $k=0$, and briefly recall the required constructions. Refer to \cite{BSSZ15, CM14} for background and proofs.
 
Observe the tautological diagram
\begin{equation}
    \xymatrix{\logDR_g(d, -\bfnu) \ar[dr]^c \ar@/_1pc/[ddr]_{Br} \ar@/^1pc/[drr]^{S}& & \\
    &\DR_g(d, -\bfnu) \ar[r]^{s} \ar[d]_{br}& \overline{\mathcal{M}}_{g,m+1}\\
    & [LM_{2g-1+m}/S^{2g-1+m}]
    }
\end{equation}
where the maps $S,s$ are source morphisms, and the maps $Br, br$ are branch morphisms. In the specific case of $k=0$ the target of the branch morphisms is a symmetric stack quotient of a Losev-Manin space. We denote by $0, \infty$ the two heavy (and distinguished) marked points in such a space. The $2g-1+m$ light points are considered unmarked. We denote by $\hat{\psi}_0\in A^1([LM_{2g-1+m}/S^{2g-1+m}])$ the $\psi$-class at the marked point $0$, by $\tilde{\psi}_0$ the $\psi$-class on $\DR_g(d, -\bfnu)$ at the unique inverse image of $0$, and by $\psi_0$ the $\psi$-class on $\overline{\mathcal{M}}_{g, m+1}$ at the corresponding mark.

We have the following facts:
\begin{itemize}
    \item by definition, the branch class in the logarithmic tautological ring of $\logDR_g(d, -\bfnu)$ is pulled back from Losev-Manin; denoting by $\Delta^c$ the sum of all boundary strata of codimension $c$ we have
    \begin{equation}\label{eq:branch}
        \branch_c = Br^\ast([\Delta^c]);
    \end{equation}
    \item since when $k = 0$ we have $Br = \pi\circ br$ , by projection formula with respect to the morphism $c$,
    \begin{equation}
       \mathsf{H}_g((d, - \bfnu), \textbf{e})_{|k = 0} := \int_{\logDR_g(d, -\bfnu)}\mathbf{\psi^e}\branch_{2g-3+|\mathbf{e}|} = \int_{\DR_g(d, -\bfnu)}\mathbf{\psi^e}br^\ast([\Delta^{2g-3+|\mathbf{e}|}]);
    \end{equation}
    \item the relation between $\hat{\psi}_0$ and $\tilde{\psi_0}$ is given by a slight generalization of Ionel lemma, which is described in this context in \cite[Proposition 2.5]{BSSZ15}, \cite[Lemma 4.2]{CM14}:
    \begin{equation} \label{eq:hattotilde}
        d\cdot \tilde{\psi}_0 = br^\ast(\hat{\psi_0});
    \end{equation}
   \item the relation between $\hat{\psi}_0$ and $\tilde{\psi_0}$, given in \cite[Lemma 2.6]{BSSZ15} or \cite[Lemma 4.3]{CM14}, simplifies because of the full ramification condition $d$:
   \begin{equation}\label{eq:tildetonot}
        \tilde{\psi}_0 = s^\ast({\psi_0});
    \end{equation}
    \item Denote by $\Delta_{r_\infty}$ the irriducible boundary divisor in $[LM_{r}/S^{r}]$ parameterizing curves where $r_\infty$ light points lie on the same component as the point $\infty$. The class $\hat{\psi}_0$ admits the following boundary expression \cite[Equation $(2)$]{BSSZ15}:
    \begin{equation}\label{eq:boulos}
        \hat{\psi}_0 = \sum_{r_\infty=1}^{r-1}\frac{r_\infty}{r} [\Delta_{r_\infty}].
    \end{equation}
\end{itemize}

Before we give a general proof of Theorem \ref{Thm:psi0_e_general_g}, we treat independently the following particular case.

\begin{lemma}\label{lam:caspart}
 \begin{equation}
    \mathsf{H}_g((d,-\bfnu), (2g-2+m, 0, \ldots, 0))|_{k=0} = \frac{d^{-1}}{2g-1+m}\cdot \mathsf{H}_g((d,-\bfnu), ( 2g-3+m, 0, \ldots, 0))|_{k=0}.
\end{equation}   
\end{lemma}
\begin{proof}
    We are in the case $r = 2g-1+m$. By dimension reasons and the natural properties of restriction of $\psi$-classes to boundary strata, we have that for $r_\infty \ge 2$
        \begin{equation}\label{eq:LMvan}
       \hat{\psi}_0^{2g-3+m}\cdot [\Delta_{r_\infty}] = 0. 
    \end{equation}
With this relation and all the ingredients developed before, we now have:
\begin{align}
&\mathsf{H}_g((d,-\bfnu), (2g-2+m, 0, \ldots, 0))|_{k=0} \nonumber\\
 =   & \int_{\DR_g(d, -\bfnu)}\psi_0^{2g-2+m}   \nonumber \\
    \stackrel{\eqref{eq:hattotilde}, \eqref{eq:tildetonot}}{=} &\int_{\DR_g(d, -\bfnu)}\frac{1}{d^{2g-2+m}}\cdot br^\ast(\hat\psi_0^{2g-2+m}) \nonumber \\
    \stackrel{\mbox{\tiny{proj.form.}}}{=} &\int_{[LM_{2g-1+m}/S^{2g-1+m}]}\frac{1}{d^{2g-2+m}}\cdot \hat\psi_0^{2g-2+m} \nonumber \\
    \stackrel{\eqref{eq:boulos}}{=} &\int_{[LM_{2g-1+m}/S^{2g-1+m}]}\frac{1}{d^{2g-2+m}}\cdot \hat\psi_0^{2g-3+m}\cdot\left(\sum_{r_\infty=1}^{2g-2+m}\frac{r_\infty}{2g-1+m}\cdot [\Delta_{r_\infty}]
    \right)
    \nonumber \\
    \stackrel{\eqref{eq:LMvan}}{=} &\int_{[LM_{2g-1+m}/S^{2g-1+m}]}\frac{1}{d^{2g-3+m}}\cdot \hat\psi_0^{2g-3+m}\cdot\frac{1}{d}\cdot \frac{1}{2g-1+m} \cdot[\Delta^1] \nonumber \\
\stackrel{\eqref{eq:hattotilde}, \eqref{eq:tildetonot}, \eqref{eq:branch}}{=} &\frac{d^{-1}}{2g-1+m} \cdot\int_{\DR_g(d, -\bfnu)}\psi_0^{2g-3+m}\branch_1 \nonumber \\
= &\frac{d^{-1}}{2g-1+m} \cdot\mathsf{H}_g((d,-\bfnu), (2g-3+m, 0, \ldots, 0))|_{k=0}.
   \end{align}
\end{proof}

\begin{remark}\label{rem:exch}
 It will be important to interpret the result of Lemma \ref{lam:caspart} in terms of weighted sums of leaky covers: the left hand side consists of a single graph with a single vertex counted with multiplicity equal to its local vertex multiplicity. The right hand side is a sum over graphs with two vertices, the second being either a trivalent rational vertex or a genus one cul de sac. Denoting by $\mathrm{v}_0, \mathrm{v}_1$ the two vertices and by $\textsc{e}$  the set of compact edges between them \footnote{Observe that $|\textsc{e}| = 1,2,3$ depending on the number of left and right edges incident to $\mathrm{v}_1$}, we have:
 \begin{equation}
  \mathsf{H}_g((d,-\bfnu), (2g-2+m, 0, \ldots, 0))|_{k=0}  = \frac{d^{-1}}{2g-1+m}\cdot \sum_{\pi} \frac{1}{|\Aut(\pi)|}\cdot \prod_{i= 0,1}\mult_{\mathrm{v}_i}\cdot \prod_{\mathrm{e}\in \textsc{e}}w(\mathrm{e}).   
 \end{equation}
\end{remark}

We are now ready to tackle the proof of Theorem \ref{Thm:psi0_e_general_g}.

\begin{proof}[Proof of Theorem \ref{Thm:psi0_e_general_g}]
We proof the theorem by induction on $e_0$, with the base case being trivially true.

Denote by $\Lambda_e$ the set of all leaky covers contributing to  $\mathsf{H}_g((d,-\bfnu), (e, 0, \ldots, 0))|_{k=0}$. For $\pi:\Gamma\to \mathbb{R} \in \Lambda_e$ denote by $\mathrm{v}_0$ the first vertex of the graph. For $\pi':\Gamma'\to \mathbb{R}\in \Lambda_{e-1}$ denote by $\mathrm{v'}_0$ the first vertex, $\mathrm{v'}_1$ the second vertex, and by $\textsc{e}'$ the set of edges between them. We denote by $\Aut(\textsc{e}')$ the set of automorphisms of the multiset $\{w(\mathrm{e'})\}_{\mathrm{e'}\in \textsc{e}'}$

For a leaky cover $\pi\in \Lambda_e$, let $$\mult_\pi=\mult_{\mathrm{v}_0}\cdot {R_\pi},$$ and for $\pi'\in \Lambda_{e-1}$, let $$\mult_{\pi'} =\prod_{i= 0,1}\mult_{\mathrm{v}_i}\cdot\frac{\prod_{\mathrm{e'}\in \textsc{e}'} w(\mathrm{e'})}{|\Aut(\textsc{e}')|} \cdot {R_{\pi'}},$$ where  $R_\star$ is a symbol defined by the two equations to denote the ``remaining'' part of the multiplicity of the leaky cover after factoring out the bit that we want to focus on.

We define a natural function 
\begin{equation}
    \varphi:\Lambda_{e-1}\to \Lambda_e
\end{equation}
which, for a graph $\Gamma'$ of a leaky cover $\pi'\in \Lambda_{e-1}$, contracts the edges in $\textsc{e}'$.

For any  fixed   $\pi:\Gamma\to \mathbb{R} \in \Lambda_e$, we make the following observations:
\begin{enumerate}
    \item For any $\pi'\in \varphi^{-1}(\pi)$, 
\begin{equation}\label{eq:istess}
    R_{\pi'} = R_\pi.
\end{equation}
\item  Let us cut $\Gamma$ after $\mathrm{v_0}$ and, for any  $\pi'\in \varphi^{-1}(\pi)$, cut $\Gamma'$ after $\mathrm{v_1}'$, and consider the fragment of graph containing $\mathrm{v_0}'$ with the multiplicity $\mult_{\pi'}/R_{\pi'}$. One can see that the weighted sum of these fragments corresponds to exchanging with boundary  divisors one of the $\psi$-classes in the integral (without any branch class)  evaluating $\mult_{\mathrm{v_0}}$ as in Remark \ref{rem:exch}. The vertex $\mathrm{v_0}$ corresponds to some genus $g_{\mathrm{v_0}}$ and some number of right ends $m_\mathrm{v_0}$, such that  the quantity $2g_\mathrm{v_0}-2+m_\mathrm{v_0} = e$. We therefore have:
\begin{equation}\label{eq:ecco}
    \mult_{\mathrm{v_0}} = \frac{d^{-1}}{e+1}\cdot \sum_{\pi'\in  \varphi^{-1}(\pi)} \frac{\mult_{\pi'}}{R_{\pi'}}.
\end{equation}
We are now ready to conclude the proof:
\begin{align}
      \mathsf{H}_g((d,-\bfnu), (e, 0, \ldots, 0))|_{k=0}  \nonumber   &= \sum_{\pi\in \Lambda_e} \mult_{\mathrm{v_0}} \cdot R_\pi \\
      \nonumber   &\stackrel{\eqref{eq:ecco}}{=} \frac{d^{-1}}{e+1}\cdot \sum_{\pi\in \Lambda_e}\sum_{\pi'\in \varphi^{-1}(\pi)} \frac{\mult_{\pi'}}{R_{\pi'}} \cdot R_\pi \\
      \nonumber   &\stackrel{\eqref{eq:istess}}{=} \frac{d^{-1}}{e+1}\cdot \sum_{\pi'\in \Lambda_{e-1}} \mult_{\pi'}\\
      &= \frac{d^{-1}}{e+1}\cdot  \mathsf{H}_g((d,-\bfnu), ( e-1, 0, \ldots, 0))|_{k=0}
      \nonumber \\
      &= \frac{d^{-e}}{(e+1)!}\cdot  \mathsf{H}_g((d,-\bfnu), ( 0, \ldots, 0))|_{k=0},
\end{align}
where the last line is obtained by plugging in the inductive hypothesis, thus concluding the proof. 
    
\end{enumerate}

\end{proof}

\section{Computer implementations} \label{Sect:Implementations}
The $k$-leaky descendants in general and some of the recursive formulas from the paper in particular have been implemented in admcycles \cite{admcycles}. To use the code, it's necessary to install the latest version from 
\begin{center}
\url{https://gitlab.com/modulispaces/admcycles}
\end{center}
directly.
The following functions are available:
\begin{itemize}
    \item \texttt{monodromy\_graphs} -- enumerates $k$-leaky tropical covers together with their associated multiplicities
    \item \texttt{k\_leaky\_descendant} -- computes $k$-leaky double Hurwitz descendants using tropical covers
    \item \texttt{J\_formula} -- provides the symbolic formulas for special cases in Theorem \ref{Thm:genus0}
    \item \texttt{H\_formula} -- calculates one-part $k$-leaky Hurwitz numbers in arbitrary genus using the recursion formula from Lemma \ref{prop-rechighergenus}
\end{itemize}

A basic example calculation is shown below:

\begin{lstlisting}
sage: from admcycles.logtaut.k_leaky_numbers import *
sage: g = 0
sage: A = [10, -2, -3, -2, -3]
sage: psi_exp = [1, 0, 0, 0, 0]  # psi^1 at the marking with weight d=10
sage: k_leaky_descendant(g, A, psi_exp)
30
\end{lstlisting}

Further example calculations can be found in \href{https://gitlab.com/modulispaces/admcycles/-/blob/master/docs/source/notebooks/Cavalieri-Markwig-Schmitt%20(2025)%20-%20One%20part%20leaky%20covers.ipynb?ref_type=heads}{this notebook}.

\bibliographystyle{alpha}
\bibliography{main}

\end{document}